\documentclass[letterpaper, 10 pt, conference]{ieeeconf}

\newtheorem{theorem}{Theorem}

\newtheorem{lemma}{Lemma}
\newtheorem{proposition}{Proposition}
\newtheorem{example}{Example}
\newtheorem{definition}{Definition}
\usepackage{tikz}
\usepackage{tikz-network}
\usetikzlibrary{shapes.geometric, arrows}
\tikzstyle{ad} = [rectangle, rounded corners, minimum width=3cm, minimum height=1cm,text centered, draw=black, fill=blue!30]
\tikzstyle{ag} = [rectangle, rounded corners, minimum width=2cm, minimum height=1cm,text centered, draw=black, fill=green!30]
\tikzstyle{cd} = [rectangle, minimum width=3cm, minimum height=1cm, text centered, draw=black, fill=orange!30]
\usepackage{xcolor}
\usepackage{color}
\input{mysymbol.sty}
\usepackage{amsmath,amsfonts}

\title{\LARGE \bf{Robust Social Welfare Maximization via Information Design in Linear-Quadratic-Gaussian Games}}

\author{Furkan Sezer and  Ceyhun Eksin
\thanks{Authors are with the Department of Industrial and Systems Engineering, Texas A\&M University, College Station, TX 77843 USA (e-mails: furkan.sezer@tamu.edu, soham.das@tamu.edu, eksinc@tamu.edu). This work was supported by NSF CCF-2008855.}}
\date{\today}
\begin{document}
\maketitle



\begin{abstract}
Information design in an incomplete information game includes a designer with the goal of influencing players' actions through signals generated from a designed probability distribution so that its objective function is optimized.  We consider a setting in which the designer has partial knowledge on agents' utilities. We address the uncertainty about players' preferences by formulating a robust information design problem against worst case payoffs. If the players have quadratic payoffs that depend on the players' actions and an unknown payoff-relevant state, and signals on the state that follow a Gaussian distribution conditional on the state realization, then the information design problem under quadratic design objectives is a semidefinite program (SDP). Specifically, we consider ellipsoid perturbations over payoff coefficients in linear-quadratic-Gaussian (LQG) games. We show that this leads to a tractable robust SDP formulation. Numerical studies are carried out to identify the relation between the perturbation levels and the optimal information structures.
\end{abstract}




\section{Introduction}

An incomplete information game is comprised of multiple agents who takes actions which maximizes their utilities depending on actions of other agents and unknown states. Incomplete information games are used to model federated edge learning \cite{miao2022}, electricity spot market \cite{verma2022}, cyber defense in EV charging \cite{zian2023} and traffic flow in communication or transportation networks \cite{brown2017studies, wu2021value}. 

\begin{figure}[h]
\centering
\begin{tikzpicture}[scale=1]
\Vertex[size=1.3,x=0.5,y=0.5]{A} 
\Vertex[size=1.3, x=2.5,y=-1.5, color=magenta]{B}
\Vertex[size=1.3, x=-2.5,y=-1.5, color=magenta]{C}
\Vertex[size=1.3,, x=3,y=1.5]{D}
\Vertex[size=1.3, x=-3,y=1.5, color=green]{E}
\Vertex[size=1.3,x=-0.7,y=-1.5, color=green]{G}
\Text[fontsize=\footnotesize,width=1.2cm,x=0.5,y=0.5]{Agent 1}
\Text[fontsize=\footnotesize,width=1.2cm,x=2.5,y=-1.5]{Agent 2} 
\Text[fontsize=\footnotesize,width=1.2cm,x=-2.5,y=-1.5]{Agent 3}
\Text[fontsize=\footnotesize,width=1.2cm,x=3,y=1.5]{Agent 4}
\Text[fontsize=\footnotesize,width=1.2cm, x=-3,y=1.5]{Agent 5}
\Text[fontsize=\footnotesize,width=1.2cm, x=-0.7 ,y=-1.5]{Agent 6}
\Vertex[size=1.8, y=3, color=yellow ]{F}
\Text[fontsize=\footnotesize,width=1.5cm,x=0, y=3]{Information Designer}
\Edge(A)(B)
\Edge(A)(C)
\Edge(A)(D)
\Edge(A)(E)
\Edge(E)(B)
\Edge(C)(G)
\Edge(G)(D)
\Edge[color=red](A)(F)
\Edge[color=red](F)(B)
\Edge[color=red](F)(C)
\Edge[color=red](F)(D)
\Edge[color=red](F)(E)
\Edge[color=red](F)(G)
\end{tikzpicture}
\caption{Information designer sends optimally designed signals on the risks of infection from an emerging infectious disease to the population with individuals who are susceptible (blue), infected (green) or recovered (magenta), so that they follow the recommended health measures, e.g. social distancing or masking that reduce the risk of an outbreak. An individual's infection or disease transmission risk is determined by its contacts (shown by black edges)--see Example 2. For instance agent 1 (susceptible) has one infected neighbor (agent 5) that it can contract the disease from.}
\label{inf_design_diagram_epidemic}
\vspace{-12pt}
\end{figure}
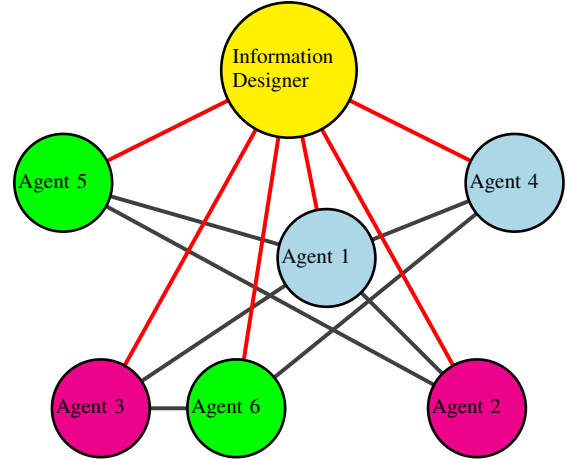

Information design problem entails decision over informativeness of signals given to agents regarding the payoff state so that induced actions maximize a system level objective. Information designer as an entity commits to an optimal probability distribution of signals conditional on payoff states before state realization (for an example in pandemic control see Fig. \ref{inf_design_diagram_epidemic}). The selected distribution maximizes the designer objective and adheres to equilibrium constraints. Various entities such as social media companies \cite{candogan2020information}, advertisements platforms \cite{emek2014} and public health agencies \cite{alizamir2020healthcare} could be considered as information designers. In control systems, information design is employed for routing games \cite{yixian2022}, Vehicle-to-Vehicle communication \cite{gould2022}, and queue management under heterogeneous users \cite{heydaribeni_anastasopoulos_2021}.

In this paper, we propose a robust optimization approach to the information design problem considering the fact that the designer cannot know the players' payoffs exactly. Indeed, while the designer may be knowledgeable about the payoff relevant random state, it may have uncertainty about the payoff coefficients of the players. For instance, in the pandemic control example (Fig. \ref{inf_design_diagram_epidemic}) above while the public health department may have near-certain information about the potential risks of a disease or intervention, it may not know how the society weights the risks and benefits in their decision-making. Here, we assume the designer has partial knowledge about players' utilities, and wants to perform information design over the payoff relevant states. 

When the payoffs of the players are unknown, the designer cannot be sure of the rational behavior under a chosen information structure. We formulate this problem as a robust optimization problem where the designer chooses the ``best'' optimal information structure for the worst possible realization of the payoffs. That is, we do not make any assumptions on the distribution of the players' payoff coefficients.  

Specifically, we assume the players have linear-quadratic payoffs with coefficients unknown by the designer. We further assume that the payoff relevant states and signals generated by the designer come from a Gaussian distribution. In this setting, we show that the robust information design with the goal to maximize social welfare can be formulated as a tractable SDP given ellipsoid  perturbations on the payoff coefficients--see Theorem \ref{theorem_model}.

In Bayesian persuasion literature, robustness is explored in worst case, online and various other settings \cite{dworczak2022preparing,hu2021robust,zu2021learning,BABICHENKO2022226,deClippel2022}. For instance, \cite{bonatti2023coordination} considers information design where the designer learns unknown utilities via auctions. Instead, here we consider the multi-player setting, i.e., information design, and assume an incomplete information game among the players. In our setting, the designer maximizes the worst-case objective under the rational behavior.

\subsection{Notation}
We use $A_{i,j}$ to denote the element in the $i$th row and $j$th column of matrix $A$. For matrices $A\in \reals^{m\times m}$ and $B\in \reals^{m\times m}$. We use $\bullet$ to represent the Frobenius product, e.g., $A\bullet B = \sum_{i=1}^{m} \sum_{j=1}^{m} A_{i,j}B_{i,j}$. We use $P^m$ and $P^{m}_{+}$ to represent the set of $m \times m$ symmetric and symmetric positive semi-definite matrices, respectively. Trace of a matrix is denoted with $\tr(\cdot)$. $I$ indicates an identity matrix. $\bbone$ is a column vector of all ones.

\section{Generic Robust Information Design Problem for Welfare Maximization}
An incomplete information game involves a set of $n$ players belonging to the set $\mathcal{N}:=\{1,\dots,n\}$, each of which selects actions $a_i \in \mathcal{A}_i$ to maximize the expectation of its individual payoff function $u_i^\theta(a,\gamma)$ where $a \equiv (a_{i})_{i \in \mathcal{N}} \in \mathcal{A} $, $\gamma \equiv (\gamma_{i})_{i \in \mathcal{N}} \in \Gamma$, and $\theta \in \Theta$ correspond to an action profile, a payoff state vector, and a payoff parameter, respectively. The payoff state of player $i$ $\gamma_i$ directly influences agent $i$'s payoff, and is unknown by the player. Agent $i$ forms expectation about the payoff state $\gamma$ based on its signal/type $\omega_i \in \Omega_i$. The payoff coefficients $\theta$ are unknown to the designer, but known to the players. 
We represent the incomplete information game given $\theta \in \Theta$ by the tuple $G_\theta:= \{\mathcal{N}, \mathcal{A}, \Gamma,\{u_i^\theta\}_{i\in \mathcal{N} }, \{\omega_i\}_{i\in \mathcal{N}}\}$. We use $\ccalG_\Theta:=\{G_\theta: \theta \in \Theta\}$ to refer to the set of games parameterized by $\theta$.

The information designer does not know that actual payoff parameter $\theta$, but knows that the game played belongs to $\ccalG_\Theta$. 
%
An information designer aims to maximize a system level objective function $f^\theta:\mathcal{A}\times \Gamma \to \reals$, e.g., social welfare, that depends on the actions of the players ($a$), and the state realization ($\gamma$) by deciding on an information structure $\zeta$ belonging to the feasible space of probability distributions on the signal space $\mathcal{Z}$ given a game with payoff coefficients $\theta$. The information structure determines the fidelity of signals $\{\omega_i\}_{i\in \mathcal{N}}$ that will be revealed to the players given a realization of the payoff state $\gamma$. 

We introduce social welfare as a design objective.
\begin{definition}[Social Welfare] \label{ex_social_welfare}
Social welfare design objective is the sum of individual utility functions, 
\begin{align}\label{eq_soc_welfare}
f^\theta (a, \gamma) &= \sum_{i=1}^{n} u_{i}^\theta(a,\gamma).
\end{align}
\end{definition}
Social welfare is a common design objective used in congestion  \cite{wu2021value}, global \cite{morris2002social}  or public goods games \cite{alizamir2020healthcare}.

A strategy of player $i$ maps each possible value of the private signal  $\omega_{i} \in \Omega_{i}$ to an action $s_{i}(\omega_{i}) \in \mathcal{A}_{i}$, i.e., $s_{i}: \Omega_{i} \rightarrow \mathcal{A}_{i}$. 
A strategy profile $s = (s_{i})_{i \in \mathcal{N}}$ is a Bayesian Nash equilibrium (BNE) with information structure $\zeta$ of the game $G_\theta$, if it satisfies the following inequality
\begin{equation}\small
E_{\zeta}[u_{i}^\theta(s_{i}(\omega_{i}), s_{-i},\gamma )|\omega_{i}] \geq  E_{\zeta}[u_{i}^\theta(a_{i}', s_{-i},\gamma )|\omega_{i}],\label{eq_class_bne}
\end{equation}
for all $ a_{i}' \in \mathcal{A}_{i}, \omega_{i} \in \Omega_{i}, i \in \mathcal{N}$,  and $s_{-i}= (s_{j}(\omega_{j}))_{j\neq i}$ is the equilibrium strategy of all the players except player $i$, and $E_{\zeta}$ is the expectation operator with respect to the distribution $\zeta$ and the prior on the payoff state $\psi$. We denote the set of BNE strategies in a game $G_\theta$ with $BNE(G_\theta).$ 

In this paper, the designer does not make any distributional assumptions on the payoff parameter $\theta$, and aims to select the best signal distribution for the worst case scenario, i.e., 
\begin{equation}\label{eq:orig-obj}
\max_{\zeta \in \mathcal{Z}} \min_{s \in BNE(G_\theta),\; \forall G_\theta \in \ccalG_\Theta} E_{\zeta}[f^\theta(s, \gamma)].
\end{equation}
Inner optimization problem in \eqref{eq:orig-obj} evaluates to the designer's objective under the worst possible payoff parameter realization and BNE actions given a signal distribution $\zeta$. The designer wants to do the best it can to maximize the system objective assuming the realization of the worst-case scenario. 

We denote the optimal solution to \eqref{eq:orig-obj} by $\zeta^*$. Given the robust optimal information structure $\zeta^*$, the information design timeline is given in the following:
\begin{enumerate}
\item Designer notifies players about $\zeta^*$
\item Realization of payoff state $\gamma$, and payoff parameter $\theta$ with subsequent draw of signals $w_{i},\, \forall i \in \mathcal{N}$ from $\zeta^*(\omega,\gamma)$ 
\item Players take action according to BNE strategies under information structure $\zeta^*$
\end{enumerate}
The generic robust information design problem in \eqref{eq:orig-obj} is not tractable in general. In the following we make assumptions on the payoff structure and the signal distribution to attain a tractable formulation. 

\subsection{Linear-Quadratic-Gaussian (LQG) Games} \label{sec_lqg}
An LQG game corresponds to an incomplete information game with quadratic payoff functions and Gaussian information structures. Specifically, each player $i \in \mathcal{N}$ decides on his action $a_{i} \in \mathcal{A}_{i} \equiv \mathbb{R}$ according to a payoff function
\begin{equation}\label{utility}
u_{i}^\theta(a , \gamma ) = - H_{i,i}a_{i}^{2} - 2 \sum_{j \neq i}H_{i,j}a_{i}a_{j} + 2\gamma_{i}a_{i} +d_{i}(a_{-i}, \gamma)
\end{equation}
where $ \mathcal{A} \equiv \mathbb{R}^{n} $ and $\Gamma \equiv \mathbb{R}^{n}$ that is a quadratic function of player $i$'s action, and is bilinear with respect to $a_i$ and $a_j$, and $a_i$ and $\gamma$. The term  $d_{i}(a_{-i},\gamma)$ is an arbitrary function of the opponents' actions $a_{-i}\equiv (a_{j})_{j \neq i} $ and payoff state $\gamma$. We collect the coefficients of the quadratic payoff function in a matrix $H = [H_{i,j}]_{n\times n}$. The payoff parameter $\theta$ unknown to the designer in \eqref{utility} is the coefficients matrix $H$.


Payoff state $\gamma$ follows a Gaussian distribution, i.e., $\gamma \sim \psi(\mu , \Sigma)$ where $\psi$ is a multivariate normal probability distribution with mean $\mu\in \mathbb{R}^n$ and covariance matrix $\Sigma$. 
Each player $i \in \mathcal{N}$ receives a private signal $\omega_{i} \in \Omega_{i} \equiv \mathbb{R}^{m_{i}}$ for some $m_i\in \mathbb{N}^+$. We define the information structure of the game  $\zeta(\omega|\gamma)$ as the conditional distribution of $ \omega\equiv (\omega_{i})_{i \in N}$ given $\gamma$.  We assume the joint distribution over the random variables $(\omega,\gamma)$ is Gaussian; thus, $\zeta$ is a Gaussian distribution. 


Next, we provide two examples of LQG games.

\begin{example}[The Beauty contest Game]
Payoff function of player $i$ is given by
\begin{equation}\label{eq:payoff_beauty}
u_{i}^\theta(a, \gamma) = -(1-\theta)(a_{i}-\gamma)^{2} - \theta(a_{i}-\bar{a}_{-i})^{2},
\end{equation}where $\theta\in [0,1]$  and $\bar{a}_{-i}= \sum_{j\neq i} a_{j}/(n-1)$ represents the average action of other players. The first term in \eqref{eq:payoff_beauty} denote the players' urge for taking actions close to the payoff state $\gamma$. The second term accounts for players' tendency towards taking actions in compliance with the rest of the population. The constant $\theta$ gauges the importance between the two terms. The payoff captures settings where the valuation of a good depends on both the performance of the company and what other players think about its value \cite{morris2002social}. 
\end{example}

\begin{example}[Social Distancing Game]
Player $i$'s action $a_{i} \in \mathbb{R}^{+} \cup \{0\}$ is its social distancing effort to avoid the infectious disease contraction/transmission (see also Fig. \ref{inf_design_diagram_epidemic}). The risk of infection depends on unknown disease specific parameters, e.g., severity, infection rate, and the social distancing actions individuals in contact with agent $i$. We define the payoff function of player $i$ as follows,
\begin{equation}\label{eq_distancing}
u_{i}^\theta(a, \gamma)=-H_{i,i} a_{i}^{2}-(1-\delta_{i}a_{i}) r_{i}(a,\gamma)
\end{equation}where the risk of infection is 
$r_{i}:=\gamma- 2\sum_{i\neq j}H_{i,j}a_{j}$,
$0<\delta_{i}<1$ is the risk reduction coefficient. In the definition of risk $r_{i}$, $\gamma$ denotes the risk rate of the disease such as infection rate or severity, and $H_{i,j}$ determines the contacts of agent $i$ and the intensity of the contacts.
First term in \eqref{eq_distancing} represents the cost of social distancing. Second term in \eqref{eq_distancing} denotes the overall risk of infection that scales with the player's social distancing efforts. 
\end{example}

Next we state the main structural assumption on perturbed LQG games. 
\begin{assumption} \label{A_structure}
We assume the following perturbation structure on the payoff matrix $H$,
\begin{equation}\label{eq_H_base_pert}
H_{i,j}=  [H_{0}]_{i,j} + v_{i,j}\epsilon_{i,j},\quad\forall i,j \in \mathcal{N}
\end{equation}where $v_{i,j} \in \mathbb{R},$ is an element of the unknown perturbation matrix $v \in \mathbb{R}^{n \times n}$ which covers a given closed and convex perturbation set $\mathcal{V}$ such that $0 \in \mathcal{V}$ and $\epsilon_{i,j}$ is the constant shift.
\end{assumption}
Assumption \ref{A_structure} means that the parameter $\theta$ in game $G_\theta$ corresponds to $H$.

\section{Robust Information Design under Finite Scenarios\label{sec_intr}}


{We will reformulate the problem in \eqref{eq:orig-obj} in order to obtain a tractable formulation. The reformulation will first entail changing the design variables from signals to actions. In order to do this, we define the distribution of actions induced by the information structure under a given strategy profile.  }
\begin{definition}[Action distribution]
An action distribution is the probability of observing an action profile $a\in \mathcal{A}$ when agents follow a strategy profile $s$ under $\zeta$, which can be computed as
\begin{equation}\label{eq_phi}
\phi(a|\gamma) = \sum_{\omega: s(\omega)=a}\zeta(\omega|\gamma). 
\end{equation}
\end{definition}
According to the definition, the probability of observing action profile $a$ is the sum of the conditional probabilities of all signal profiles $\omega$ under $\zeta$ that induce action profile $a$ given the strategy profile $s$.

We denote the set of {\it equilibrium} action distributions induced by BNE strategies under an information structure $\zeta \in \ccalZ$ for game $G_\theta$ as 
\begin{equation} \label{eq_eq_action_dist}\small
C(\ccalZ) = \{ \phi :\phi \textrm{ satisfies \eqref{eq_phi} for }  s\in BNE(G_\theta) \textrm{ given } \zeta \in \ccalZ \}.
\end{equation}

The designer can recommend actions instead of sending signals to each player, if the designer knew the payoff coefficient $\theta$. In such a case, the players would follow the recommended actions because they would satisfy the obedience condition as per the revelation principle, see \cite[Proposition 1]{Bergemann2019}. However, this principle does not apply in the setting where $\theta$ is adversarially chosen. To overcome this issue, we assume the obedience condition is only satisfied in the worst case scenario. We detail our approach first in the finite-scenario case, where $\theta$ can take finite set of values. 

%


We begin by stating the BNE condition in \eqref{eq_class_bne} by a set of linear constraints for LQG games given the payoff coefficients $H$. 

\begin{lemma} \label{lem_linear}
Define the covariance matrix $X \in P^{+}_{2n}$ as follows:
\begin{equation}\label{eq_x_defn}
X:=\begin{bmatrix}
var(a) & cov(a, \gamma)\\
cov(\gamma, a) & var(\gamma)
\end{bmatrix}.  
\end{equation}
For a given payoff matrix $H$, the BNE condition in \eqref{eq_class_bne} can be written as the following set of equality constraints, 
\begin{align} 
\sum_{j \in \mathcal{N}}H_{i,j} X_{i,j}-X_{i,n+i} =0 , \quad i\in \mathcal{N}\label{eq_affected_inequality_app}
\end{align}
where $X_{i,j}=cov(a_i,a_j)$ for $i\leq n$, and $j\leq n$, and $X_{i,n+i}=cov(a_i,\gamma_i)$. 
\end{lemma}
\begin{proof}
See Appendix.    
\end{proof}
The condition in \eqref{eq_affected_inequality_app} ensures that $X$ is a Bayesian correlated equilibrium (BCE), see \cite{Bergemann2019} for a definition.

In the following, we express the robust information design problem under a finite set of scenarios as a mixed integer SDP. 
\begin{proposition}[Finite-case] \label{prop_finite}
Let the design objective $f^\theta(a,\gamma)$ be quadratic in its arguments with the coefficients stored in matrix $F\in \reals^{2n \times 2n}$, i.e., $f^\theta(a,\gamma) = [a\; \gamma]^T F [a\; \gamma]$. Suppose Assumption \ref{A_structure} holds, and assume the design objective coefficients do not depend on $H$. Consider a finite perturbation vector with $C$ scenarios, and let $v_{c} \in \reals^{n \times n}$ refer to perturbation vectors corresponding to one of the scenarios $c\in \mathcal{C}= \{1,\dots,C\}$. We can express the robust information design problem in \eqref{eq:orig-obj} as the following mixed-integer SDP:
\begin{align}
&\min_{y_{c} \in \{0,1\},\forall c \in \{1,2,..,C\}}  \max_{X \in P^{2n}_{+} }  \; F\bullet X \label{eq_semi_first_K}\\
\text{s.t.} \; &y_{c} (R_{0,l}\bullet X+ \sum_{(i,j) \in \mathcal{Y}_{l}} [v_{c}]_{i,j}\epsilon_{i,j} X_{i,j})= 0,  \nonumber\\ &\forall l \in \mathcal{N}, c\in \mathcal{C}\label{eq_uncertain_general_form_semi_K}\\
&\sum_{c=1}^{C}y_{c} =1, \\
&M_{k,l}\bullet X = cov(\gamma_{k},\gamma_{l}), \quad\forall k,l \in \mathcal{N} \text{ with } k\leq l, \label{eq:rmodend_semi_K} 
\end{align} 
where $X$ is defined in \eqref{eq_x_defn}, $R_{0, l}= [[R_{0, l}]_{i,j}]_{2n\times 2n} \in P^{2n}, l \in \mathcal{N}$
is given as: 
\begin{equation}\label{eq_R_epsilon}
[R_{0, l}]_{i,j}=\begin{cases}
[H_{0}]_{l,l}&if \quad  i = j = l,  \\
[H_{0}]_{l,j}/2  & if \quad  i = l, 1\leq j \leq n, j\neq l, \\
-1/2 &if \quad i = l, j= n + l,\\
[H_{0}]_{i,l}/2 &if\quad  j = l, 1\leq i \leq n, i\neq l \\
-1/2 &if \quad j = l, i= n + l,\\
0 & \text{otherwise,}
\end{cases}
\end{equation}
 $M_{k,l}= [[M_{k,l}]_{i,j}]_{2nx2n} \in P^{2n},  k \in \mathcal{N}$ is given as: 
\begin{equation}\label{eq_M_matrix}
[M_{k,l}]_{i,j}=\begin{cases}
1/2 \quad\text{ if } k<l, i=n+k, j=n+l\\
1/2 \quad\text{ if } k<l, i=n+l, j=n+k\\
1 \quad\text{ if } k=l, i=n+k, j=n+l\\
0 \quad \text{otherwise,}
\end{cases} \qquad 
\end{equation}
and $[v_{c}]_{i,j}$ refer to the elements of the perturbation vector with 
\begin{align} \label{eq_Y_set}
\mathcal{Y}_{l}&:=\{\{i,j\}: i = j = l \lor i = l, 1\leq j \leq n, j\neq l \nonumber\\& \lor j = l, 1\leq i \leq n, i\neq l  \}.
\end{align}

\end{proposition}
\begin{proof}
We can express the expected objective using the Frobenius product as follows,
\begin{align}
E_{\phi}[f(a,\gamma)] &= E_\phi\big[ \begin{bmatrix}
a^{T}, &\gamma^{T}
\end{bmatrix} F \begin{bmatrix}
a \\ \gamma
\end{bmatrix} \big ], \\
&= F \bullet X
\end{align}where 
$
F=\begin{bmatrix}
[F]_{1,1} & [F]_{1,2}\\ [F]_{1,2} & [F]_{2,2}
\end{bmatrix} \in P^{2n},
$ and note that $[F]_{i,j}$ denotes the $i,j$th $n\times n$ submatrix.

Let $c^*$ be the worst-case scenario from the perspective of the designer. The designer chooses $X^*$ that maximizes its objective $F \bullet X$ subject to rational behavior of players in the worst case scenario. As per Lemma \ref{lem_linear}, we have 
\begin{align} 
& \sum_{j \in \mathcal{N}}H_{i,j} X^*_{i,j}-X^*_{i,n+i} =0,  \quad \forall i\in \mathcal{N} \\
&\sum_{j \in \mathcal{N}}([H_0]_{i,j} +[v_{c^*}]_{i,j}\epsilon_{i,j})  X^*_{i,j}-X^*_{i,n+i} =0,   \forall i\in \mathcal{N}. \label{eq_bce_k}
\end{align}
We rewrite \eqref{eq_bce_k} in terms of matrices $R_{0,l}, \forall l\in \mathcal{N}$ as in \eqref{eq_R_epsilon} and $X$ as in \eqref{eq_x_defn} to obtain \eqref{eq_uncertain_general_form_semi_K}. Minimization over $y_{c}, \{1,2,..,C\}$ enforces the constraint $c^*$ among the set of constraints in  \eqref{eq_uncertain_general_form_semi_K} to be selected.   
Constraint \eqref{eq:rmodend_semi_K} corresponds to the assignment of $var(\gamma)$ to $[X]_{2,2}.$ Constraint \eqref{eq:rmodend_semi_K} is not affected by perturbations to $H.$ 
\end{proof}

According to the formulation in \eqref{eq_semi_first_K}-\eqref{eq:rmodend_semi_K}, the solution can entail finding the covariance matrix $X$ that maximizes $F\bullet X$ for each scenario $c=1,\dots,C$, and then picking the smallest among them. We note that an alternative equivalent formulation can entail $C$ covariance matrices, i.e., $X_1,\dots,X_{C}$, and leave out the integer variables $\{y_{c}\}_{c=1,\dots,C}$. 

We use the scenario-based formulation \eqref{eq_semi_first_K}-\eqref{eq:rmodend_semi_K} to motivate the tractable robust design formulations under ellipsoid and interval formulations. For illustration purposes, consider $C=2$ scenarios. Assume scenario $c=1$ is the worst case scenario, i.e., $y_1=1$ and $y_2=0$. In such a case, $X^*$ will satisfy the BNE condition \eqref{eq_bce_k} for $c=1$ exactly while the BCE condition will be approximately satisfied for $c=2$. Specifically, we have
\begin{align}
& \sum_{j \in \mathcal{N}}([H_{0}]_{i,j} +[v_{2}]_{i,j}\epsilon_{i,j})  X^*_{i,j}-X^*_{i,n+i}= \\&     \sum_{j \in \mathcal{N}}([H_{0}]_{i,j} +[v_{2}]_{i,j}\epsilon_{i,j}+[v_{1}]_{i,j}\epsilon_{i,j}-[v_{1}]_{i,j}\epsilon_{i,j})  X^*_{i,j}-X^*_{i,n+i} \\
&= \sum_{j \in \mathcal{N}} ([v_{2}]_{i,j}\epsilon_{i,j} - [v_{1}]_{i,j}\epsilon_{i,j})  X^*_{i,j} >0
\end{align}
We can interpret this relation as the optimal solution to \eqref{eq_semi_first_K}-\eqref{eq:rmodend_semi_K} $X^*$ being induced by an approximate BNE for the good scenario $c=2$. That is, $X^*$ is not necessarily incentive compatible with players' realized payoffs. 
In the following, we leverage this observation to develop robust convex program for social welfare objective when the perturbation set is an ellipsoid.  

\section{An SDP Formulation for Social Welfare Maximization via Information Design\label{sec_trac}}




Under convex uncertainty sets, the number of scenarios $C$ goes to infinity. Thus, we cannot enforce exact BCE explicitly for the worst-case scenario, and annul the other cases using integer variables as is done in \eqref{eq_semi_first_K}-\eqref{eq:rmodend_semi_K}. Instead, we relax the BCE constraint in \eqref{eq_affected_inequality_app} as follows
\begin{align}
\sum_{j \in \mathcal{N}}H_{i,j} X_{i,j}-X_{i,n+i} \leq \alpha, \quad i\in \mathcal{N}.\label{eq_alpha_affected}
\end{align}
where $\alpha>0$ is a finite large-enough constant. Consider the following ellipsoid uncertainty subsets $\mathcal{V}_{l} \subset \mathcal{V}, \forall l \in \mathcal{N}$: 
\begin{equation}\label{eq_ellipsoidal_perturbation}
\mathcal{V}_{l} = \text{Ball}_{\rho} = \{v_{l} \in \mathbb{R}^{2n-1}: \lvert \lvert v_{l} \rvert \rvert_{2} \leq \rho  \}, \quad \forall  l \in \mathcal{N}.
\end{equation}

We take the social welfare (Example \ref{ex_social_welfare}) as the designer's objective $f^\theta(a,\gamma)$, which depends on the payoff matrices $\theta \equiv H$. 

\begin{theorem}\label{theorem_model}
Assume $H$ is given by \eqref{eq_H_base_pert} and perturbation vectors $v_{l}, 
\forall l \in \mathcal{N}$ exhibit ellipsoid uncertainty \eqref{eq_ellipsoidal_perturbation} and the objective is social welfare maximization with 
\begin{equation}\label{eq_soc_coef}
F= \begin{bmatrix}
-H & I\\I & O
\end{bmatrix}\text{ and } F_{0}= \begin{bmatrix}
-H_{0} & I\\I & O
\end{bmatrix}.
\end{equation}
The robust convex program under the welfare maximization objective is as follows:
\begin{align}
&\max_{X \in P^{2n}_{+},t} \quad t \label{model_robust_sdp_tract_welfare}\\ 
\text{s.t. }& F_{0}\bullet X - \frac{n^{2}\rho}{2n-1} \sqrt{ \sum_{i=1}^{n} \sum_{j=1}^{n} (\epsilon_{i,j} X_{i,j})^{2}} \geq t, \label{eq_soc_tractable}\\
&R_{0,l}\bullet X+\rho \sqrt{ \sum_{(i,j) \in \mathcal{Y}_{l}} (\epsilon_{i,j} X_{i,j})^{2}}  \leq \alpha, \quad \forall l \in \mathcal{N} \label{eq:rmod2}\\
&M_{k,l}\bullet X = cov(\gamma_{k},\gamma_{l}), \quad\forall k,l \in \mathcal{N} \text{ with } k\leq l. \label{model_robust_sdp_tract_welfare_end}
\end{align}
 where matrices $R_{0, l}$ and $M_{k,l}$ are as defined in \eqref{eq_R_epsilon} and \eqref{eq_M_matrix}, respectively.
\end{theorem}
\begin{proof}
See \cite{sezer2021social} on how to express the social welfare objective in \eqref{ex_social_welfare} using \eqref{eq_soc_coef}, and in form $F\bullet X$.
We start by writing the social welfare objective constraint $F\bullet X \geq t$ under ellipsoid uncertainty:
\begin{equation}\label{eq_social_welfare_semi_inf}
F\bullet X = F_{0}\bullet X -\sum_{i=1}^{n} \sum_{j=1}^{n} v_{i,j}\epsilon_{i,j} X_{i,j}  \geq t
\end{equation}
Here we consider all elements of perturbation matrix $v$ for ellipsoid perturbations:
\begin{equation}\label{eq_ellipsoidal_perturbation_F}
\mathcal{V} = \text{Ball}_{\rho} = \{v \in \mathbb{R}^{n\times n}: \lvert \lvert v \rvert \rvert_{2} \leq \frac{n^{2}\rho}{2n-1}  \}.
\end{equation}
Using \eqref{eq_ellipsoidal_perturbation_F}, we will obtain a robust counterpart for semi-infinite constraint \eqref{eq_social_welfare_semi_inf}. We start with writing \eqref{eq_social_welfare_semi_inf} as a perturbation maximization problem:
\begin{equation}\label{eq_social_welfare_pert_max}
\max_{\lvert \lvert v\rvert \rvert \leq \rho}  \sum_{i=1}^{n} \sum_{j=1}^{n} v_{i,j}\epsilon_{i,j} X_{i,j} \leq  F_{0}\bullet X-t
\end{equation}
Solution to \eqref{eq_social_welfare_pert_max} is the tractable robust constraint \eqref{eq_soc_tractable}. 
    
Next, we substitute $H$ with \eqref{eq_H_base_pert} into \eqref{eq_alpha_affected}:
\begin{align}\small
\sum_{j \in \mathcal{N}}([H_{0}]_{ij} + v_{ij}\epsilon_{ij}) cov(a_{i},a_{j})-cov(a_{i},\gamma_{i}) &\leq \alpha, \forall i\in \mathcal{N}\label{eq_affected_inequality_sub}
\end{align}
We can rewrite \eqref{eq_affected_inequality_sub} in terms of matrices $R_{0,l}, \forall l\in \mathcal{N}$ and $X$ as in \eqref{eq_x_defn}:
\begin{align}
R_{0,l}\bullet X+\sum_{(i,j) \in \mathcal{Y}_{l}} v_{i,j}\epsilon_{i,j}X_{i,j}&\leq \alpha , \quad \forall l\in \mathcal{N}.\label{eq_affected_inequality_X}
\end{align}

Similarly, we write the perturbation maximization problem over uncertain constraint \eqref{eq_affected_inequality_X} under ellipsoid uncertainty as 
\begin{align}
\max_{\lvert \lvert v_{l} \rvert \rvert \leq \rho}\sum_{(i,j) \in \mathcal{Y}_{l}} v_{i,j}\epsilon_{i,j}X_{i,j}&\leq \alpha-R_{0,l}\bullet X , \quad \forall l\in \mathcal{N}\label{eq_affected_inequality_prob}
\end{align}
where $\ccalY_l$ is given by \eqref{eq_Y_set}. 
Solution to \eqref{eq_affected_inequality_prob} give us the tractable constraint \eqref{eq:rmod2}.

Constraint \eqref{model_robust_sdp_tract_welfare_end} enforces assignment of known covariance matrix of payoff states, $cov(\gamma)$ to the respective place in $X$. 
\end{proof}
The equilibrium constraints given in \eqref{eq:rmod2} make sure the recommended action distribution is an approximate BCE for every realization of the payoff coefficients matrix.

\section{Numerical Experiments\label{sec_nume}}
We consider a designer that wants to maximize the social welfare of $n=5$ players. The designer knows the perturbed payoff coefficients given as follows: $[H_{0}]_{i,i}=5$ for $i\in\{1,\dots, 5\}$, and $[H_{0}]_{i,j}=-1$ for $i\neq j$, $i,j \in \{1,2,..,5\}$. 
%
The variance of the unknown payoff state $\gamma$ is given as follows: $var(\gamma)_{i,i}=5$ for $i=\{1,\dots, 5\}$, and $var(\gamma)_{i,j}=0.5$ for $i \neq j, \;i,j \in \{1,2,.,5\}$.
%
%
We consider ellipsoid perturbations with $\rho\in \{0.7,1,1.3,..,3.4\}$ and let $\alpha =0.1$. Given the setup, we solve the robust convex program \eqref{model_robust_sdp_tract_welfare}-\eqref{model_robust_sdp_tract_welfare_end} in order to obtain the robust optimal information design $X^*$.

We analyze the effects of shifts $\epsilon_{i,j}$ defined in \eqref{eq_H_base_pert} by assuming the diagonal elements and off-diagonal elements of shift matrix are homogeneous, i.e.,  $\epsilon_{i,i}=\epsilon_1$ and $\epsilon_{i,j} = \epsilon_2$ for all $i,j=1,\dots,n$ for constants $\epsilon_1$ and $\epsilon_2$. 

\begin{figure}[ht]
\begin{tabular}{c}
\includegraphics[width=0.4\textwidth]{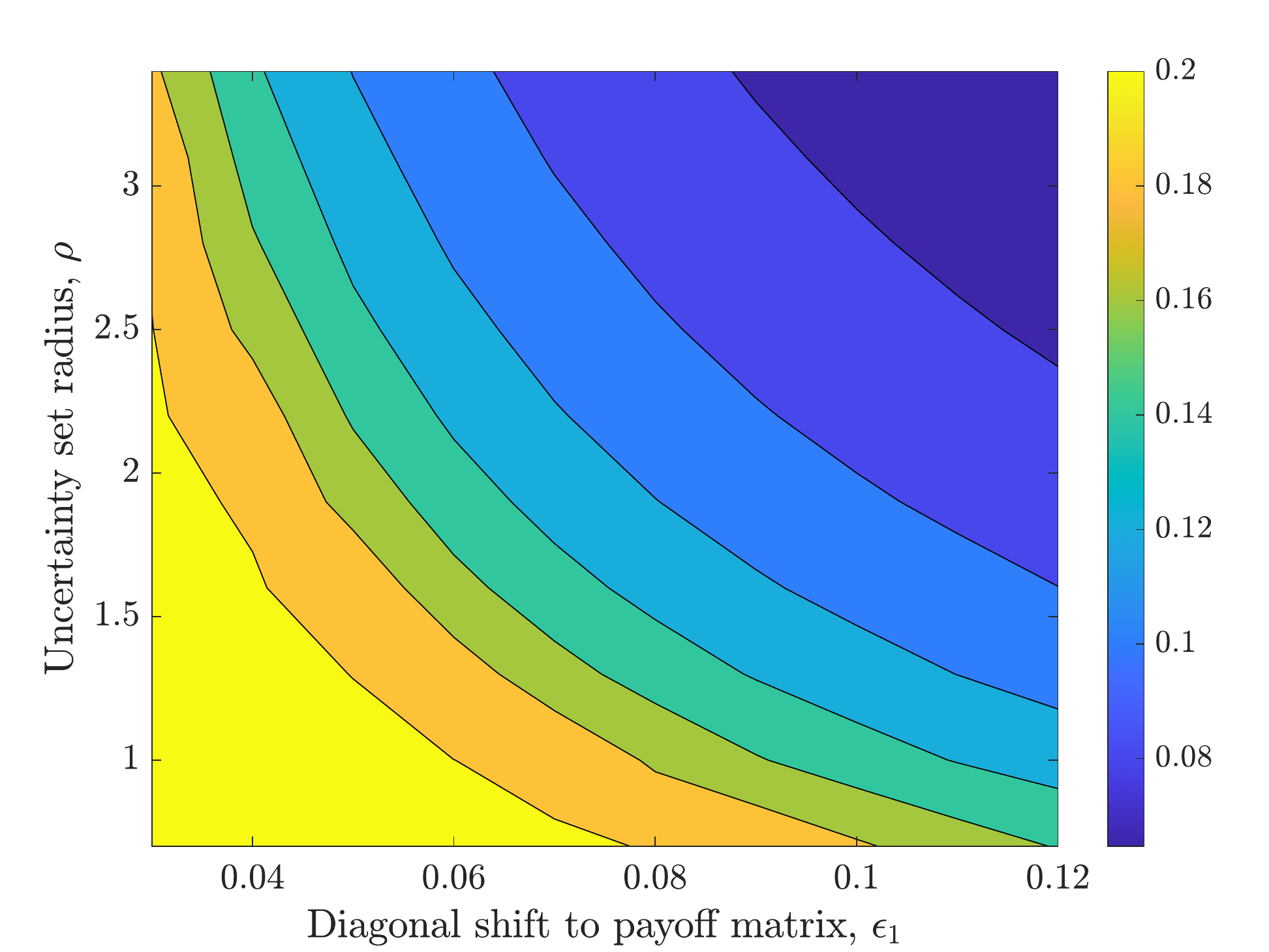} \\
\small (a) $\left \vert \left \vert  X^* - X_{no} \right \vert \right \vert_{F}$ \\
\includegraphics[width=0.4\textwidth]{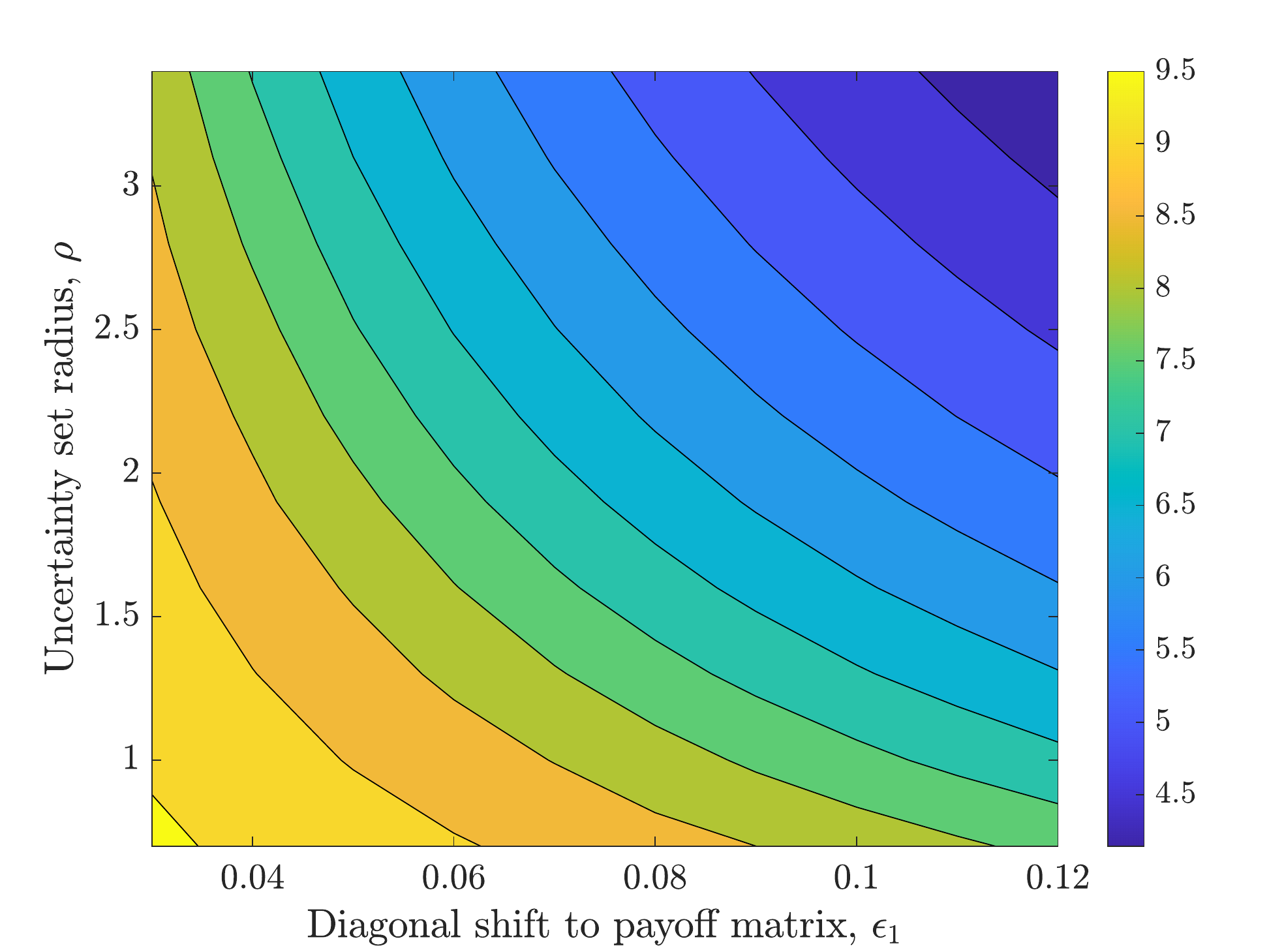} \\ 
\small (b) Optimal objective value 
\end{tabular}
\caption{\small Contour plots of (a) normalized Frobenius matrix norm distance  $\left \vert \left \vert  X^* - X_{no} \right \vert \right \vert_{F}$ between the optimal covariance matrix ($X^*$) and no information disclosure covariance matrix ($X_{no}$), and (b) optimal objective value with respect to uncertainty ball radius $\rho$ and diagonal shift $\epsilon_{1}$ to coefficient matrix $H$ under a symmetric supermodular game with social welfare objective. Optimal solution $X^*$ approaches to no information disclosure as $\rho$ and $\epsilon_{1}$ increase.}
\label{fig_diag}
\end{figure}

In order to systematically analyze the effects of the shifts, we fix the off-diagonal shifts to a small value $\epsilon_2=0.001$, and vary the diagonal shift $\epsilon_{1} \in \{0.03,0.04,0.05,..,0.12 \}$. Fig. \ref{fig_diag}(a) shows that as the uncertainty ball radius $\rho$ and diagonal shift $\epsilon_{1}$ increases, the optimal information structure remains a partial information disclosure but gets closer to the no information disclosure. Fig. \ref{fig_diag}(b) shows that social welfare decreases under increasing uncertainty.

We can discuss Fig. \ref{fig_diag}  in terms of the beauty contest game, which is a supermodular game. If we consider the common goods in the beauty contest game as a stock, we see that a social welfare maximizing information designer, i.e. the company whose stock is traded releases less information about stock price $\gamma$, when the uncertainty about its shareholder's payoff coefficients $H$ increases.

\section{Conclusion}
The paper considered the problem of designing information structures in incomplete information games when the designer does not know the game payoffs exactly. This is a common situation in many real-world settings, where the game payoffs are often uncertain due to various factors such as incomplete information, imperfect modeling, or unknown parameters. Specifically, we considered information design for the setting when the unknown payoff parameters are adversarially chosen. For the robust information design problem, we developed a tractable SDP formulation given quadratic payoffs, Gaussian signal distributions, ellipsoid perturbations to the unknown payoff parameters, and social welfare as the design objective. Numerical experiments show that the designer would choose to reveal less information about the payoff states to the players as its uncertainty about the players' payoffs grow. This suggests that in situations where the game payoffs are highly uncertain, it may be more optimal to not disclose any information at all rather than risk providing misleading information.

\bibliographystyle{ieeetr}
\bibliography{ref}

\appendix

\section{Proof of Lemma \ref{lem_linear}}
We start with writing the first order condition equivalent to \eqref{eq_class_bne} for a given $\theta \equiv H$: 
\begin{align}\label{eq_first_order}
&E_{\zeta} \bigg[\frac{\partial}{\partial a_{i}} u_{i}^\theta (s(\omega), \gamma) | \omega_{i} \bigg] = -2 H_{i,i} s_{i}(\omega_{i})\nonumber \\&-2 \sum_{i \neq j} H_{i,j} E_{\zeta}[s_{j}| \omega_{i}] + 2 E_{\zeta}[\gamma_{i}| \omega_{i}] =0
\end{align}
We solve \eqref{eq_first_order} for the best response $s_{i}(\omega_{i}), \forall i \in \mathcal{N}$:
\begin{equation}\label{eq_strategy}
 H_{i,i}s_{i}(\omega_{i}) = \sum_{i \neq j} H_{i,j} E_{\zeta}[s_{j}| \omega_{i}] - E_{\zeta}[\gamma_{i}| \omega_{i}], \quad i\in\mathcal{N} 
\end{equation}
%
We look for an equilibrium strategy of the form given below:
\begin{equation}\label{eq_linear_BNE_app}
s_{i}(\omega_{i}) = \bar{a}_{i} +b_{i}^{T} (\omega_{i} - E_{\zeta}[\omega_{i}]) , \quad \forall i \in \mathcal{N}, 
\end{equation}
where $\bar{a}_{i}$ and $b_{i}^{T}, \forall i \in \mathcal{N}$ are constants and constant vectors, respectively.  We plug \eqref{eq_linear_BNE_app} into the first order condition \eqref{eq_strategy}:
\begin{align}
&\sum_{j \in \mathcal{N}} H_{i,j} E[\bar{a}_{j} +b_{j}^{T} (\omega_{j} - E_{\zeta}[\omega_{j}]) | \omega_{i}=\bar{\omega}_{i}] = 	E[\gamma_{i}| \omega_{i}=\bar{\omega}_{i}],
\end{align}
$\forall \bar{\omega}_{i} \in \Omega_{i},  i \in \mathcal{N}$. Via conditional expectation rule over multivariate normal distribution, we obtain following:
\begin{align}
&\sum_{j \in \mathcal{N}} H_{i,j} (b_{j}^{T} cov(\omega_{j}, \omega_{i}) var(\omega_{i})^{-1} (\bar{\omega}_{i}-E_{\zeta}[\omega_{i}] ) + \bar{a}_{j}) \nonumber\\ &=  E[\gamma_{i}] +cov(\omega_{i},\gamma_{i})^{T} var(\omega_{j})^{-1} (\bar{\omega}_{i}- E_{\zeta}[\omega_{i}]), \label{eq_general_equlibrium}
\end{align}
$\forall \bar{\omega}_{i} \in \Omega_{i}, i \in \mathcal{N}$.
Vectors $b_{i}, i \in \mathcal{N} $ and constants $\bar{a}_{i}, i \in \mathcal{N}$ are determined by following set equations when we separate \eqref{eq_general_equlibrium} into respective parts:
\begin{align}\label{eq_BCE_detailed}
\sum_{j \in \mathcal{N}} H_{i,j} b_{j}^{T} cov(\omega_{j}, \omega_{i}) var(\omega_{i})^{-1} &=  cov(\omega_{i},\gamma_{i})^{T} var(\omega_{j})^{-1}, \nonumber\\ &   \forall \bar{\omega}_{i} \in \Omega_{i}, i \in \mathcal{N} 
\end{align}
\begin{equation}
\sum_{j \in \mathcal{N}} H_{i,j} \bar{a}_{j} =  E[\gamma_{i}], \quad i \in \mathcal{N}.\label{eq_constant_prior}
\end{equation}

We divide both sides of \eqref{eq_BCE_detailed} by $var(\omega_{i})^{-1} $ and obtain the following set of equations: 
\begin{align}\label{eq_bne_b}
\sum_{j \in \mathcal{N}} H_{i,j} b_{j}^{T} cov(\omega_{j}, \omega_{i})  =  cov(\omega_{i},\gamma_{i}), \quad i \in \mathcal{N} 
\end{align}
For scalar signals $\omega_i\in \reals$, if we let $b_i=1$ and $\bar a_i = E_\zeta [\omega_i]$ for $i\in \ccalN$, then we have $a_i =\omega_i$ by \eqref{eq_linear_BNE_app}. Moreover, the set of equations in \eqref{eq_bne_b} is equivalent to \eqref{eq_affected_inequality_app}. 

\end{document}